\documentclass[preprint,12pt]{revtex4-2}

\usepackage{graphicx} % Required for inserting images
\usepackage{color}
\usepackage{xcolor}%
\usepackage{amssymb,amsmath,environ}
\usepackage{empheq}
\usepackage{url}
\usepackage{booktabs}

\numberwithin{equation}{section}

\def\CTcomment[#1]{\textcolor{teal}{#1}}
\def\MCTcomment[#1]{\textcolor{blue}{#1}}
\def\GLcomment[#1]{\textcolor{red}{#1}}
\def\AScomment[#1]{\textcolor{violet}{#1}}

\def\bu{\mathbf u}
\def\bw{\mathbf w}

\def\bk{\mathbf k}

\newcommand{\Gr}{\mathcal{G}}

\begin{document}

\title{Spreading of pathological proteins through brain networks: a case study for Alzheimer's disease}

\author{Germana Landi}
\affiliation{Department of Mathematics, University of Bologna}
\email[Corresponding author:]{ germana.landi@unibo.it}
\author{Arianna Scaravelli}
\affiliation{Department of Mathematics, University of Bologna}
\author{Maria Carla Tesi}
\affiliation{Department of Mathematics, University of Bologna}
\author{Claudia Testa}
\affiliation{Department of Physics and Astronomy ``Augusto Righi'', University of Bologna}
%%\title{Spreading of pathological proteins through brain networks: a case study for Alzheimer's disease}
%\author{Germana Landi, Arianna Scaravelli, Maria Carla Tesi}
%\address{Department of Mathematics, University of Bologna}
%\author{Claudia Testa
%}
%\address{Department of Physics and Astronomy ``Augusto Righi'', University of Bologna
%}

% corresponding author
%\corraddr{germana.landi@unibo.it}

\begin{abstract}
Mathematical modeling offers a valuable approach to understanding Alzheimer's disease (AD) given its complexity, unknown causes, and lack of effective treatments.
Models, once validated, offer a powerful tool to test medical hypotheses that are otherwise difficult to verify directly. Our focus here is on elucidating the spread of misfolded $\tau$ protein, a critical hallmark of AD alongside A$\beta$ protein, taking also into account the synergistic interaction between the two proteins. We consider distinct modelling choices, all employing network frameworks for protein evolution, differentiated by their network architecture and diffusion operators. By carefully comparing these models against clinical $\tau$ concentration data, gathered through advanced multimodal analysis techniques, we show that certain models replicate better the protein's dynamics. This investigation underscores a crucial insight: in modeling complex pathologies, the precision with which the mathematical framework is chosen is crucial, especially when validation against clinical data is considered decisive.


\end{abstract}

\maketitle

%\AScomment[Arianna]
%\MCTcomment[Maria Carla]
%\GLcomment[Germana]
%\CTcomment[Claudia]
%
%Rivista:
%\url{https://www.aimspress.com/mbe/article/6838/special-articles}
%
%deadline: spostata al 31 ottobre 2025

%\input{introduzione}

\section{Introduction}
Alzheimer's disease (AD) is an irreversible and incurable neurodegenerative disorder, progressively eroding memory and cognitive function in over 50 million people worldwide, a number projected to surge in the coming years (World Alzheimer Report 2025) \cite{ADreport}. While AD's exact causes remain elusive, two proteins, amyloid-beta ($A\beta$) and tau ($\tau$), are central to its emergence and development. Though naturally present in the brain, in AD they abnormally aggregate into $A\beta$ plaques and neurofibrillary tangles (NFTs), respectively, which are hallmarks of the disease \cite{ittner2011amyloid,bloom2014amyloid,tatarnikova2015beta,bennett2017enhanced,small2008linking}.
These pathological proteins don't spread uniformly; instead, they follow distinct spatiotemporal patterns \cite{ossenkoppele2018discriminative,ahmed2014novel}. NFTs typically emerge in the entorhinal cortex before spreading throughout the brain, while $A\beta$ plaques initially form in temporal and frontal areas \cite{braak1991neuropathological,grothe2017vivo}. Recent mathematical research highlights the crucial interplay between $A\beta$ and $\tau$, emphasizing its consideration for effective therapies \cite{busche2020synergy,thompson2020protein,bertsch2021macroscopic}. For an updated review on mathematical models we refer the interested reader to \cite{moravveji2024scoping,carbonell2018mathematical}.

When it comes to modelling the spreading of proteins in the brain, a very convenient setting is represented by networks. Indeed, different brain regions are ``connected" from a structural point of view by bundles of fibers, and this structure can be very naturally modelled by means of brain networks.
A quite commonly used network nowdays is represented by a connectome \cite{mcnab2013human,bertsch2023,
torok2023connectome,fornari2019prion,raj2012network}. From a mathematical point of view, a connectome is a weighetd graph, in which the vertices represents parcellated regions of gray matter formed by clusters of neurons, which
share similarities in cytoarchitecture, functional activity, and structural connections to other regions, while
edges represent the connectivity between the regions, which can differ accordingly to the weights. In \cite{bianchi2024network}  
we have been using
the Budapest Reference Connectome v3.0 \cite{szalkai2017parameterizable} as a basic connectome:
it is a weighted graph that has been obtained by averaging  the tractography of diffusion tensor images of 477 healthy subjects of the Human Connectome Project \cite{mcnab2013human}.

We stress that the vertices of such a graph do not have associated coordinates that determine their position in space, since it is precisely an averaged graph. This will have consequences when we talk about “distances”, which must be understood in an intrinsic sense and not in the classical Euclidean metric sense.
We use this \emph{structural connectome} to build new connectomes with the aim of best reproducing some key features of the phenomena we were interested in describing.
More precisely, in \cite{bianchi2024network} we introduced two ``new" connectomes for the spreading of the two proteins $A\beta$ and $\tau$. Since $A\beta$ protein is known to travel on short distances in the extracellular space \cite{meyer2008rapid}, we constructed what we called the ``intrinsic proximity graph", as follows. Given two vertices on the graph, we connect them only if they are close in an intrinsic sense, i.e. only if there already exists an edge connecting them, shorter that a chosen 
threshold. There is no ``vicinity" in an Euclidean metric sense, since 
the vertices of the graph do not have associated spatial coordinates. We then assign weights to the edges, that are stronger for intrinsically close vertices.
We believe that this modelling choice is reliable, so also in this work
 the $A\beta$ protein will spread along this graph, by means of the graph Laplacian operator associated to it.

Unlike $A\beta$, $\tau$ is believed to travel over long distances, in the intracellular space \cite{dujardin2020tau,goedert2017propagation}. To take into account this feature, we introduced in \cite{bianchi2024network}
what we called the ``cumulative connectome".
This connectome is constructed by connecting two vertices of the structural connectome if there is at least one path between them that is shorter than a fixed  threshold. The weight of the edge between the two vertices is determined by the sum of the lengths of all paths connecting them. This aggregation of path lengths is what gives the connectome its 'cumulative' nature.
In \cite{bianchi2024network}, we explored several approaches to model $\tau$ spreading and compared the resulting $\tau$ concentration values with clinical data. The comparison focused on significant brain regions, where "significant" refers to regions we have identified through statistical analysis of clinical data, highlighting differences in $\tau$ aggregation between normal and AD brains.
%In \cite{bianchi2024network}
%we considered several different ways to model $\tau$ spreading, and we compared the results on $\tau$ concentration with clinical data related to brain regions which exhibits a significant difference in $\tau$ aggregration between normal and AD brains.  
%in areas of the brain considered relevant when it comes to the AD pathology.  
Usually, diffusion is modelled via a Laplacian (in our case a graph Laplacian) or via a convolution operator, especially when long distances are concerned. While for the diffusion of $A\beta$ protein we simply used the graph Laplacian operator associated to the intrinsic proximity graph, for the spreading of $\tau$ we compared different modelling choices. 
Unsurprisingly, the results in \cite{bianchi2024network} showed that 
the choice of both the graph on which diffusion takes place and of 
the mathematical operator used to model it is crucial in achieving a good match with clinical data.
The most satisfactory modelling choice turned out to be the one based on the graph Laplacian associated with the cumulative connectome.
This suggests that the cumulative nature of the cumulative connectome is best suited to describe $\tau$ diffusion, reflecting, in some sense, the “intrinsic” geometry of the brain.
%, when it comes to studying the spread of the protein through the graph Laplacian.
In the wake of this consideration, in this paper we consider a convolution operator with kernel built using cumulative paths between brain regions.
This choice proves to be decisive to obtain results improved with respect to the ones obtained in \cite{bianchi2024network}, where a different choice for the kernel was made. In addition to the introduction of a new kernel, in this paper we have
increased the number of significant regions obtaining satisfactory results.

Summarizing, the main purpose of this work is to show that both the topology of the graph over which the $\tau$ protein spreads and the operator chosen for diffusion are crucial ingredients for a mathematical model producing a good match with clinical data. The topology of the graph captures the intrinsic geometry of the phenomenon described, and therefore, must also be taken into account in the construction of the operator modelling the spreading, which somehow becomes an operator “intrinsic” to the geometry itself. 
%In our study, we see that the cumulative connectome,also used for building an intrinsic kernel for the convolution operator, is the one best reproducing the dynamics of $\tau$ spreading.
%
The paper is organized as follows: in Section 2 we describe the mathematical model we study, presenting the equations in a quite general form, that is without explicitly specifying the form of the operators adopted for diffusion or of the graph considered, but explaining carefully all the terms appearing in the equations. In Section 3 we describe the different connectomes (i.e. weighted graphs) we use, and we explicit the form of the different diffusion operators we consider. In Section 4 the procedures we have adopted to obtain reliable clinical data are explained. 
In Section 5, we present the results obtained by numerical simulations. By comparisons with clinical data we show the crucial role of the cumulative kernel, and the improvement we obtain here with respect to paper \cite{bianchi2024network}. Finally, Section 6 concludes the paper with a discussion and directions for future research.

\section{The mathematical model}

We are interested in modeling the dynamics of the two proteins $A\beta$ and $\tau$, mainly in toxic conformations at least for what concerns $\tau$, on a macroscopic scale, that is when the whole brain is considered. In particular, we would like to understand to what extent the (possibly different) strength of the connections between different brain regions could influence the spread of the proteins in the brain, and more specifically in brain networks.

The physical concept of brain network finds its ideal mathematical correspondence in the concept of graph. In this setting, the strength of the connections in the network is represented by the weights associated to the edges of the graph and the structure of the network itself by the connectivity matrix associated with the graph. In Alzheimer's disease, $\tau$ and $A\beta$ proteins have different biological mechanisms of propagation: $A\beta$ can spread over short distances across brain tissue \cite{meyer2008rapid} while $\tau$ propagates via a prion-like, neuron-to-neuron mechanism \cite{}. Moreover, their spreading follows different patterns:  $\tau$ primarily spreads through neuronal connections and $A\beta$ initially accumulates in specific cortical regions \cite{}. Therefore, in our model we assume that the world in which proteins live, that is, interact and travel, is represented by two different graphs, one for $\tau$ and one for $A\beta$. These graphs have the same vertices, corresponding to specific regions of the brain, but different edges, corresponding to different connections between regions.

                  The full system of equations we are going to study is the one considered  in 
\cite{bianchi2024network}, to which we refer for a detailed description, given by (all the constants, i.e. greek letters, are positive and $\odot$ denotes the Hadamard product between vectors):
\begin{subequations} \label{eq:system}
 \begin{empheq}[left=\empheqlbrace]
 {align}
  & \epsilon \dfrac{d \bu_1(t)}{d t} =  - \gamma_1 L_{A\beta} \bu_1(t) +  C_{\bu_1}
 -  \alpha \bu_1(t) \odot \sum_{j=1}^{3} \bu_j(t) -\sigma_1 \bu_1(t) \label{eq:u1}
  \\
 & \epsilon \dfrac{d \bu_2(t)}{d t} =  - \gamma_2 L_{A\beta}  \bu_2(t)
 +\dfrac{\alpha}{2} \bu_1(t) \odot \bu_1(t)
 -  \alpha \bu_2(t)\odot\sum_{j=1}^{3} \bu_j(t)-\sigma_2 \bu_2(t) \label{eq:u2} \\
 & \epsilon \dfrac{d \bu_3 (t)}{d t} = \frac{\alpha}{2}\sum_{3\leq j+k < 6}  \bu_j(t) \odot \bu_k(t) -\sigma_3 \bu_3(t) \label{eq:u3} \\
 & \dfrac{d \bw(t)}{d t} =  \gamma_3 K[ \bw ] +
C_{\bw} \big( \bu_2(t) -U_{\bw} \big)^+ + s_{\bw}(t) - \sigma_4 \bw(t) \label{eq:w}
\end{empheq}
\end{subequations}
The vectors $\bu_1, \bu_2$ and $\bu_3$ represent the molar concentrations of $A\beta$ monomers, dimers and plaques, respectively, in each node.
The $\epsilon$ parameter in front of their evolution equations stands for a fast dynamics \cite{meyer2008rapid}. Monomers and dimers of amyloid-beta do diffuse on short distances in the extra-cellular space, while plaques don't,
and this is modelled by the graph Laplacian operator $L_{A\beta}$, which is not present in the equation for plaques. All the terms with $\sigma_i$, $i=1,\ldots,3$  represent clearance, while the remaining terms model aggregation phenomena and $C_{\bu_1}$ is a source term for $A\beta$ monomers. Notice that the factor $\frac{\alpha}{2}$ in equations (\ref{eq:u2}) and (\ref{eq:u3}) is there for statistical reasons, to avoid double counting the same term.
In the first three equations, relative to the $A\beta$ protein, there is no interaction with $\tau$. The last equation describes the evolution, with a slow dynamics, of toxic $\tau$ in each node represented by the vector $\bw$. The last two terms in the equation model a time dependent source  and a clearance (with coefficient $\sigma_4$), respectively. As we will see later, the source of toxic $\tau$ protein will be localized in the enthorinal cortex \cite{braak2011alzheimer}. The second term in (\ref{eq:w}) encapsulates the interaction between the two proteins: a concentration of toxic $A\beta$ above a certain threshold  acts as a source for $\tau$. 
The first term in (\ref{eq:w}) is crucial in this paper, as it was in \cite{bianchi2024network}. In a possible prion-like mechanism of $\tau$ spreading, the misfolded $\tau$ protein acts as a template for the healthy one, which in turn misfolds \cite{busche2020synergy}. These phenomena occur essentially in the intracellular space, along axon bundles, over long distances between non-adjacent brain regions. 
In order to understand the temporal pattern of $\tau$ spreading, we consider in this paper two possible  spreading mechanisms on graph: diffusion and convolution.

When the spreading of $\tau$ is modeled via diffusion, chosen an appropriate graph, we have
\begin{equation}\label{eq:diff}
    K[ \bw ] = L_\tau \bw
\end{equation}
where $L_\tau$ is the usual graph Laplacian defined via the adjacency matrix \cite{grigor2018introduction}. 

When convolution is used to model $\tau$ spreading, as it is often customary for interactions over long distances, the $K$ operator is defined as follows:
\begin{equation}\label{eq:conv}
    K[ \bw ] = U \big( \hat{\bk}_\tau \odot \hat{\bw} \big), \quad \hat{\bk}_\tau=U^* \bk_\tau, \quad  \hat{\bw}=U^* \bw
\end{equation}
where $\bk_\tau$ is a kernel defined on the chosen graph, $U$ is a matrix whose columns are the eigenvectors of the graph Laplacian 
matrix and $U^*$ is the adjoint matrix.
We recall that the operator $U$ is used to define the graph Fourier transform \cite{ricaud2019fourier}; therefore, $\hat{\bk}_\tau$ and $\hat{\bw}$ are the graph Fourier transform of $\bk_\tau$ and $\bw$ and the convolution operator \eqref{eq:conv} is defined as the inverse graph Fourier transform of the element-wise product between $\hat{\bk}_\tau$ and $\hat{\bw}$.
We refer the reader to \cite{bianchi2024network,shuman2013emerging} for a deeper discussion on the graph convolution operator. 

The results in \cite{bianchi2024network} show that the choice of the graph on which the mathematical operators \eqref{eq:diff} and \eqref{eq:conv} are defined is crucial to appropriately reproduce the clinical data.  A key contribution of this work consists in having devised a particular graph and defined an appropriate convolution kernel through which we have shown a better matching with clinical data with respect to more standard and commonly used graphs.

\section{Proteins spreading through brain's connectomes}\label{sect:connectome}
As the spreading mechanisms of the $\tau$ and $A\beta$ proteins are different, they require distinct structures to be accurately described. Therefore, as already mentioned in Section 2, we use two different graphs: one for
equations \eqref{eq:u1}-\eqref{eq:u3} and one for equation \eqref{eq:w}, both based on the human ``connectome”. 
%\MCTcomment[The ‘human connectome’ is the complete map of neural connections and pathways between neurons in the human brain. The Human %Connectome Project (HCP) is a major scientific initiative that aims to map this network of anatomical and functional connections in the brain.]

The human connectome \cite{}, a map of brain structural networks, is crucial for understanding brain diseases. Representing the brain as a graph allows us to assess information processing and transfer. This representation effectively explains two fundamental brain properties: integration and segregation, which balance specialized processing with global coordination. In the graph model, brain regions are vertices connected by edges representing biological connections. 
The weights of the edges
represent the intensity of the connection between brain regions. The Budapest Reference Connectome
v3.0 [4] is a parameterizable consensus brain graph widely used in the literature. It has been
derived from the connectomes of 477 people, each computed from MRI datasets of the Human Connectome
Project. At the website \url{https://pitgroup.org/connectome/} a high-resolution version with
1015 vertices can be downloaded; the connections between brain regions can be weighted by different
factors such as the ratio between the number and the length of fiber tracts linking the vertices or the length of these fiber tracts (the choice of the preferred weight is up to the user).
We will refer to this connectome as \emph{structural connectome} since it has been reconstructed using  diffusion tensor images. 

Mathematically, the structural connectome is a weighted graph $\Gr=(V,E)$ with $N=1015$ vertices and weights given by
\begin{equation}\label{eq:ws}
  w_\textsc{nl}(i,j) = \frac{n_{ij}}{\ell_{ij}}
\end{equation} or  by
\begin{equation}\label{eq:ws2}
  w_\textsc{l}(i,j) = \ell_{ij} ,
\end{equation}
where $n_{ij}$ is the mean number of fibers connecting vertices $i$ and $j$, and $\ell_{ij}$ is the mean lenght of such fibers, that we will call \textit{fiber lenght}. We denote by $\Gr_\textsc{l}$ and $\Gr_\textsc{nl}$ the structural connectomes with edge's weights 
$w_\textsc{l}$, and  
$w_\textsc{nl}$, respectively. Both graphs share the same set of vertices and edges, corresponding to those of the original connectome 
$\Gr$; only the weights of the edges differ.
%
%The \emph{Intrinsic Proximity Connectome} and \emph{Cumulative Connectome}, introduced in \cite{bianchi2024network},  are built by using the structural connectome.
%
We define the convolution operator on the structural connectomes  $\Gr_\textsc{nl}$ and $\Gr_\textsc{l}$ by using their respective  Laplacians to construct the operator $U$, which defines the graph Fourier transform. 
We emphasize that, in the following, distances are to be intended in an intrinsic sense, i.e. referring to fiber lengths, not in a 
geometric sense, i.e. referring to the euclidean distance between two vertices. This choice is the only reasonable one, since we do not know the exact position in space of the vertices of the graph, being the graph obtained by averaging data from different subjects.

Starting from the structural connectome, in \cite{bianchi2024network} we introduced two new connectomes: the intrinsic proximity connectome and the cumulative connectome.
The \emph{intrinsic proximity connectome} 
%as been originally introduced in \cite{bianchi2024network} as 
is a weighted graph $\Gr_\textsc{p}$ obtained by connecting two nodes $i$ and $j$ of $\Gr$ 
if the fiber length $\ell_{ij}$ between them is sufficiently small. For such a graph, we fix a threshold $R_\textsc{p}$ and we define the weights as follows:
\begin{equation*}
  w_\textsc{p}(i,j) = \left\{
                        \begin{array}{ll}
                          e^{-\ell_{ij}^2/\delta_\textsc{p}} & \hbox{if $\ell_{ij}\leq R_\textsc{p}$} \\
                          0 & \hbox{otherwise}
                        \end{array} ,
                      \right.
\end{equation*}
where $\delta_\textsc{p}\in\mathbb{R}^+$ is a fixed parameter. In this way, we connect two vertices only if they are close in an intrinsic sense, instead of a ``geometric" proximity measured with some Euclidean distance, and we assign stronger weights to
intrinsically closer vertices. 
Since the $A\beta$ protein is known to spread over short distances in the brain \cite{meyer2008rapid}, we define equations \eqref{eq:u1}-\eqref{eq:u3} on the intrinsic proximity connectome and take its Laplacian as operator $L_{A\beta}$. 

The \emph{cumulative connectome}, also introduced in \cite{bianchi2024network}, is a weighted graph $\Gr_\textsc{c}$ with the same vertices as $\Gr$. We recall that, in a weighted graph, the length of a path is obtained as the sum of the weights of all edges along the path. Therefore, in the sequel we will denote as $\ell_\textsc{l}(p)$ and $\ell_\textsc{nl}(p)$ the length of a path $p$ in $\Gr_\textsc{l}$ and $\Gr_\textsc{nl}$, respectively.
%$\Gr$ measured as the sum of the weights $w_\textsc{l}$ and $w_\textsc{nl}$, respectively. 
Given a positive parameter $R_\textsc{c}$, two vertices $i$ and $j$ are directly connected in $\Gr_\textsc{c}$ if in $\Gr_\textsc{l}$ there is at least one path $p$ between $i$ and $j$ whose length $\ell_\textsc{l}(p)$ is at most $R_\textsc{c}$.
%less or equal $\delta$. 
We call such a path an \emph{admissible path}. The weight $w_\textsc{c}(i,j)$ of the edge between $i$ and $j$  
is given by the sum of the lengths $\ell_\textsc{nl}$ of all the admissible paths between $i$ and $j$:
\begin{equation*}
  w_\textsc{c}(i,j) =    \sum_{p \in \mathcal{A}_{ij}}\ell_\textsc{nl}(p) , \quad 
  \mathcal{A}_{ij}=\left\{p \; | \;  p  \text{ is a path in $\Gr_\textsc{l}$ from $i$ to $j$ and } \ell_{\textsc{l}}(p) \leq R_\textsc{c} \right\} .
\end{equation*}
We refer the reader to \cite{bianchi2024network} for a deeper description of the cumulative connectome. The cumulative connectome turned out to be the fundamental tool for modeling through diffusion the spreading of the $\tau$ protein, which is known to move along long distances in the brain. In fact, an edge between two vertices exists in $\Gr_\textsc{c}$ if the corresponding brain regions are joined by axonal paths with length controlled by the parameter $R_\textsc{c}$.
The weights of the connections depend on how many and how long the axonal fibers are.
For this reason, in this paper we use the cumulative connectome as the graph to define the operator $K$ when it models the spreading via diffusion 
%as in \eqref{eq:diff}. 
and we choose for $L_\tau$ in \eqref{eq:diff} the graph Laplacian of the cumulative connectome \cite{bianchi2024network}.

The convolution operator \eqref{eq:conv} needs selecting a proper kernel $\bk_\tau$. In \cite{bianchi2024network}, a kernel constructed by using the structural connectome with weights $w_\textsc{l}$ was considered, giving unsatisfactory results compared to clinical data. In this work, we introduce a new kernel that takes into account the very nature of the cumulative connectome. The cumulative connectome encodes \emph{all} long-range anatomical connections between brain regions and the strenght of these connections depends on the number and length of the axonal fibers. 
This approach captures the overall strength of the structural connectivity between brain regions. The cumulative nature of this connectome, which integrates all available pathways rather than relying solely on direct connections, represents a novel feature of the model and has shown promising results in aligning model predictions with clinical data \cite{bianchi2024network}.
Building on the concept of cumulative connectivity, this work introduces a \emph{cumulative kernel} that incorporates the full spectrum of connections between brain regions. This cumulative kernel, capturing not only direct links but also indirect pathways, provides a richer and more integrative representation of inter-regional interactions. In this way, it offers a powerful tool for modeling the long-range spreading of $\tau$ protein across the brain through graph convolution.
%
%We use this information to construct the convolution kernel as a signal defined over the graph nodes, capturing the cumulative connectivity profile of each region.
Importantly, the convolution is performed in the structural connectomes $\Gr_\textsc{l}$ and $\Gr_\textsc{nl}$, using the eigenvectors of their graph Laplacians as the basis for the graph Fourier transform.
This formulation allows us 
to use the spectral properties of the structural connectome to model biologically informed spatial patterns and to encode long-range propagation dynamics of $\tau$ protein into the kernel.
%This formulation allows us to encode biologically informed spatial patterns into the kernel and to use the spectral properties of the structural connectome to model the long-range propagation dynamics of $\tau$ protein.
%
%
Let $i$ be a vertex in $\Gr_\textsc{l}$; we denote by $\mathcal{M}(i)$ the set of vertices that are connected to $i$ through an admissible path, i.e.,
\begin{equation*}
    \mathcal{M}(i) = \left\{ j \in \Gr_\textsc{l} \; \middle| \; \exists \text{ admissible path } p \text{ from } i \text{ to } j \right\}.
\end{equation*}
We first define a vector $\mathbf{d} \in \mathbb{R}^N$, where each component $d_i$ quantifies the cumulative connectivity of node $i$ as follows:
\begin{equation*}
    d_i = \sum_{j \in \mathcal{M}(i)} \sum_{p \in \mathcal{A}_{ij}} \ell_\textsc{l}(p), \quad 
    \mathcal{A}_{ij} = \left\{ p \; \middle| \; p \text{ is a path in } \Gr_\textsc{l} \text{ from } i \text{ to } j \text{ and } \ell_\textsc{l}(p) \leq R_\textsc{c} \right\}
\end{equation*}
Based on this cumulative measure, we define the cumulative kernel $\bk_\tau \in \mathbb{R}^N$ as a node-dependent Gaussian-like weighting function:
\begin{equation}\label{eq:kernel}
    k_\tau(i) = e^{-d_i^2 / \delta_\textsc{k}}, \quad i = 1, \ldots, N
\end{equation}
where $\delta_\textsc{k} \in \mathbb{R}^+$ is a fixed scaling parameter.  In this formulation, nodes with large $d_i$ are assigned lower kernel values, whereas
nodes with small $d_i$  receive an higher weight.

\section{Materials and methods}\label{sect:methods}
We aim to compare the model output with clinical data in brain regions that have undergone significant changes due to Alzheimer's disease.
We use clinical imaging data to identify such brain regions and to estimate $\tau$ protein concentration values within them. 
The clinical imaging data for $\tau$ 
%used to assess %the model's %reliability
were obtained from the ADNI Initiative (ADNI) database (\url{https://adni.loni.usc.edu/}) 
%as described in \cite{bianchi2024network}
. 
In this paragraph, we briefly recall how the data were obtained, but we refer the reader to \cite{bianchi2024network} for a more detailed description.  Our study included 261 participants from ADNI3: 238 cognitively normal (CN) subjects and 23 with Alzheimer’s disease (AD). All participants were selected based on the availability of both PET scans using the [\textsuperscript{18}F]-AV1451 tracer and corresponding T1-weighted MR images acquired with a 3T Siemens scanner (MPRAGE sequence), to reduce variability related to scanner type and radiopharmaceutical. PET and MRI acquisitions were required to be within 3 months of each other to minimize the impact of disease progression between the two acquisitions.

The MRI structural images were pre-processed, segmented and parcellated with FreeSurfer 6.0 (\url{http://surfer.nmr.mgh.harvard.edu/}), in order to subdivide the brain volumes into a set of 83 anatomical cortical and sub-cortical region of interest (ROIs).
Since the nodes of the brain graphs are in correspondence with the 83 ROIs, we were able to compare $\tau$ concentrations relative to the 83 ROIs in the PET images with those estimated by the mathematical models.
To assess differences in [\textsuperscript{18}F]-AV1451 uptake by the brain regions between AD and CN subjects, PET images were co-registered to each subject’s corresponding structural MRI using PETSurfer (\url{https://surfer.nmr.mgh.harvard.edu/fswiki/PetSurfer}). 
Then, $\tau$ concentration was computed in each ROI. For absolute quantification, regional $\tau$ values were normalized to the cerebellum. Finally, a robust statistical approach, including normality assessment, non-parametric testing, and correction for multiple comparisons, was used to identify ROIs with significantly different $\tau$ concentrations between CN and AD groups. %
Our statistical analysis revealed that the distribution of \textsuperscript{18}F-AV1451 uptake is significantly different between AD and CN groups in 29 ROIs.
Table \ref{tab:signif} lists the ten most significant ROIs along with their average $\tau$ concentrations. For each ROI, the table also indicates the corresponding functional network to which it belongs. 
\begin{table}[h]
\caption{$p$-value and $\tau$ concentration (mean $\pm${ st.dev.}) in AD and CN for the most significant ten ROIs.\label{tab:signif}}
\begin{center}
\begin{tabular}{llccc}
\hline
ROI           & Network               & significance        & AD $\tau$ conc. & CN $\tau$ conc. \\ \hline
Fusiform region       &   Occipital    & $8.3 \cdot10^{-9}$  & 1.6 $\pm$ {0.3}\   & 1.2 $\pm$ {0.1} \\
Inferior temporal region   & Temporal & $1.6\cdot10^{-8}$   & 1.7 $\pm$ {0.5}\   & 1.2 $\pm$ {0.2} \\
Middle temporal region   &  Temporal  & $2.9 \cdot10^{-6}$  & 1.6 $\pm$ {0.5}\  & 1.2 $\pm$ {0.2} \\
Lingual region        &   Occipital    & $3.1 \cdot10^{-6}$  & 1.4 $\pm$ {0.3}\   & 1.1 $\pm$ {0.1} \\
Lateral orbitofrontal region & Limbic & $5.4 \cdot 10^{-6}$ & 1.5 $\pm$ {0.4}   & 1.2 $\pm$ {0.2} \\
Amygdala                 &  Limbic  & $6.5\cdot10^{-6}$   & 1.4 $\pm$ {0.4}   &   1.2 $\pm$ {0.1} \\ 
Temporalpole region   &    Temporal   & $1.3\cdot10^{-5}$   &  1.4 $\pm$ {0.3}          &  1.1 $\pm$ {0.2} \\  
Enthorinal region      &  Limbic     & $2.1\cdot10^{-5}$   &  1.4 $\pm$ {0.4}               &    1.1 $\pm$ {0.2}               \\
Parsorbitalis region      & Frontal  & $2.6\cdot10^{-5}$   &    1.5 $\pm$ {0.5}                &     1.2 $\pm$ {0.3}               \\
Lateraloccipital region   & Occipital  & $2.7\cdot10^{-5}$   &     1.5 $\pm$ {0.6}              &   1.2 $\pm$ {0.2}               \\
\hline
\end{tabular}
\end{center}
\end{table}

\section{Results}\label{results}
A series of numerical simulations was carried out to determine whether the diffusion-based or convolution-based formulation of the operator $K$ provides a better representation of the clinical imaging data.
All simulations were conducted using Matlab R2021a on an Intel Core i5 processor with
2.50 GHZ and a Windows operating system. The codes used for the current experiments can be made available
upon reasonable request to the authors. 
\subsection{Experimental setting}
We used the same experimental setting described in \cite{bianchi2024network}.
The differential system \eqref{eq:system} was numerically solved using Matlab’s \texttt{ode45}, which implements the Dormand–Prince Runge–Kutta method with adaptive time stepping \cite{dormand1980family}. %
The system was integrated from the initial time $t_0=0$ to the final time $t_f=500$ chosen to ensure stabilisation of the solutions.
The source term for misfolded $\tau$ was set in the enthorinal region which is the brain region exhibiting the earliest $\tau$ deposits in Alzheimer's disease; i.e,
\begin{equation}
    (s_\bw(t))_j = \left\{
  \begin{array}{ll}
    1 & j\in\hbox{ entorhinal region} \\
    0 & \hbox{otherwise}
  \end{array}
\right.
\end{equation}
where $(s_\bw(t))_j$ denotes the $j$th component of $s_\bw(t)$. 
%
% --------------------------- VALORI DEI PARAMETRI
%
In our numerical simulations, we have fixed the values of all parameters of system \eqref{eq:system} as reported in Table \ref{tab:parametri}, except $\gamma_3$ and $C_{\bw}$. These values have been chosen, on the basis of different tests, as the ones giving the best approach to asymptotic equilibrium in a pathological situation. 
\begin{table}[h]
    \caption{Fixed parameter values for system \eqref{eq:system}.\label{tab:parametri}}
\begin{center}
   \begin{tabular}{llllllllll}
\toprule
$\gamma_1$ & $\gamma_2$ & $\alpha$ & $C_{\bu_1}$ & $U_{\bw}$ & $\sigma_1$ & $\sigma_2$ & $\sigma_3$ & $\sigma_4$\\
\midrule
0.001 & 0.001 & 0.1 & 0.05  & 0.01 & 0.1 & 0.1 & 0.1 & 0.11 \\
\bottomrule
\end{tabular} 
\end{center}
\end{table}
%
%
% \begin{table}[H]
%     \caption{Fixed parameter values for system \eqref{eq:system}.\label{tab:parametri}}
% \begin{center}
%    \begin{tabular}{lllllllllll}
% \toprule
% $\gamma_1$ & $\gamma_2$ & $\alpha$ & $C_{\bu_1}$ & $C_{\bw}$ & $U_{\bw}$ & $\sigma_1$ & $\sigma_2$ & $\sigma_3$ & $\sigma_4$\\
% \midrule
% 0.001 & 0.001 & 0.12 & 0.1 & 0.5 & & 0.1 & 0.1 & 0.1 & 0.11 \\
% \bottomrule
% \end{tabular} 
% \end{center}
% \end{table}
%
%
As described in Section \ref{sect:connectome}, we model the spreading of the $A\beta$ protein on the intrinsic proximity connectome $\Gr_\textsc{p}$, where the values $R_\textsc{p}=25$ and $\delta_\textsc{p}=1.5 \cdot 10^2$ have been fixed. When the spreading of the $\tau$ protein is modeled through a diffusion operator, the underlying network structure is represented by the cumulative connectome $\Gr_\textsc{c}$ with $R_\textsc{c}=30$. Finally, when employing a convolution operator to model the spreading of $\tau$ protein, we use the structural connectome $\Gr$ and compare the two different weighting schemes: one based on fiber lengths, and the other on the ratio between the number of fibers and their lengths. This leads us to define two distinct Fourier transform operators: $U_\textsc{nl}$, corresponding to the structural graph with weights $w_\textsc{nl}$, and $U_\textsc{l}$, corresponding to the graph 
with weights $w_\textsc{l}$. These operators encode different spectral properties of the connectome and are used to analyze the $\tau$ protein dynamics under convolution. In both cases, the kernel $\bk_\tau$ is defined by \eqref{eq:kernel} with $\delta_\textsc{k}=10^{-4}$.
\subsection{Clinical Deterioration Pattern and Model Evaluation}

In this work, we adopt the same methodology used in our previous paper \cite{bianchi2024network} to evaluate the agreement between model predictions and clinical data. This approach is based on comparing the \textit{clinical deterioration pattern}, derived from $\tau$ concentrations in significative brain regions (as identified in Section 4), with the pattern predicted by the mathematical model. 
While in our previous study we considered only six ROIs, selected according to a statistical threshold with $p$-values up to the order of $10^{-6}$, in the present work we extend the analysis to ten ROIs, including regions with $p$-values up to the order of $10^{-5}$. The ten ROIs considered are listed in Table~\ref{tab:signif}, along with the functional networks to which they belong.
Moreover, to assess the models' ability to distinguish affected from unaffected regions, we also considered the $\tau$ concentration in the sensorimotor network, including the paracentral, postcentral, precentral and superior frontal regions, which is typically unaffected in AD pathology. 

\begin{figure}[h!]
    \centering
\includegraphics[width=0.35\linewidth]{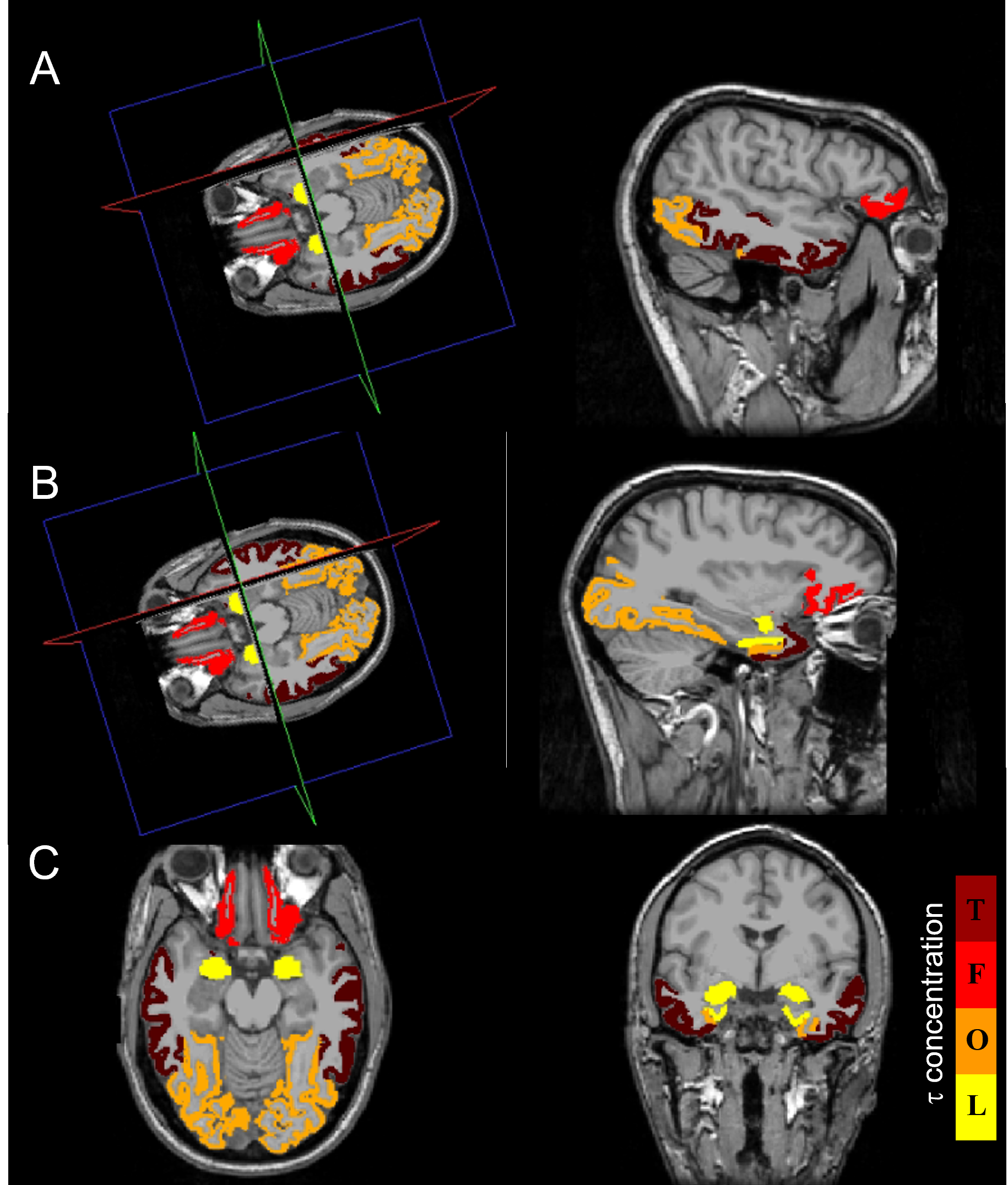}
    \caption{Mean $\tau$  concentration values $ w^{(*)}_T,\, w^{(*)}_O,\, w^{(*)}_L,\, w^{(*)}_F $ in the temporal, occipital, fusiform and limbic networks as listed in Table \ref{tab:signif}. The colorbar reflects the ordering of $\tau$ values, from high (dark red) to low (yellow). A: 3D view of the brain  and the corresponding sagittal plane. B: 3D view of the brain  and the corresponding sagittal plane. C: corresponding axial and coronal view.
 }
    \label{fig:brain}
\end{figure}

We denote by $ w^{(*)}_T,\, w^{(*)}_O,\, w^{(*)}_L,\, w^{(*)}_F $ the mean $\tau$ values averaged over the significant ROIs belonging respectively to the temporal, occipital, limbic and frontal networks as detailed in Table \ref{tab:signif}. That is, each $ w^{(*)}_X $ represents the mean $ \tau $ value across all significant ROIs within the corresponding network $ X $. 
Figure \ref{fig:brain} represents
the mean $\tau$ concentration values $ w^{(*)}_T,\, w^{(*)}_O,\, w^{(*)}_L,\, w^{(*)}_F $ in the four networks.
In addition, we denote by $w^{(*)}_S$ the mean $\tau$ value in the sensorimotor netwok.
%
%The average $\tau$ values across the significative ROIs grouped in the temporal, occipital, limbic and frontal networks, as detailed in Table \ref{tab:signif}, and across the sensorimotor network are denoted respectively by \( w^{(*)}_T, w^{(*)}_O, w^{(*)}_L, w^{(*)}_F, w^{(*)}_S \), and 
The values $ w^{(*)}_T,\, w^{(*)}_O,\, w^{(*)}_L,\, w^{(*)}_F, w^{(*)}_S $ form a reference profile derived from clinical data and ordered from highest to lowest $\tau$ concentration, they yield a string that we define as the \textit{clinical deterioration pattern}. This pattern reflects the spatial progression of neurodegeneration in AD, with higher $\tau$ levels indicating more severe deterioration. In \cite{bianchi2024network} we considered six ROIs (grouped into three networks) obtaining the clinical deterioration pattern represented by the string TLOS since \( w^{(*)}_T > w^{(*)}_L > w^{(*)}_O > w^{(*)}_S \). In this work, with ten ROIs ( grouped into four networks), the corresponding string becomes TFOLS. 
We notice that when ten ROIs are considered not only a new network (the Frontal one) occurs but also additional regions contribute to the values \( w^{(*)}_T, w^{(*)}_O, w^{(*)}_L, w^{(*)}_S \).

To evaluate the ability of the mathematical model \eqref{eq:system}  to reproduce the clinical deterioration pattern, we numerically solve the differential system and compute the average $\tau$ concentrations over the same brain networks considered in the clinical analysis. This is possible since in the graph (the connectome) we know precisely to which region belongs each vertex. These values are then used to construct a string, based on the descending order of $\tau$ concentrations.

We then compare the model’s deterioration pattern to the clinical reference one using the \textit{Hamming distance (HD)}, which counts the number of mismatches between the two strings. A lower HD indicates a closer match to the clinical progression of the disease. This approach provides a quantitative framework to assess how well the mathematical model captures the spatial dynamics of $\tau$ pathology in AD, both in affected and preserved brain regions.
\subsection{Numerical Results}
Here, we present the numerical results obtained by solving the differential system \eqref{eq:system} for the following choices of the operator $K$: diffusion on the cumulative connectome $\Gr_\textsc{c}$ and convolution on the structural connectomes $\Gr_\textsc{l}$ and  $\Gr_\textsc{nl}$.
Tables~\ref{tab:numerical_results_tlos} and~\ref{tab:numerical_results_tfols} show the string corresponding to the minimal Hamming distances obtained for each configuration, together with  the parameter values that allow the model to reproduce it. Table~\ref{tab:numerical_results_tlos} refers to the case with six ROIs and the deterioration string TLOS, while Table~\ref{tab:numerical_results_tfols} refers to the case with ten ROIs and the string TFOLS.
The parameters $C_{\mathbf{w}} $ and $\gamma_3$ were heuristically identified in each case to minimize the Hamming distance between the clinical and simulated deterioration patterns. For each parameter, we observed intervals of values that yield the same minimal Hamming distance; in the tables, we report one representative value per parameter.
\begin{table}[h!]\ 
\centering
\caption{Deterioration patterns and model parameters for different choices of the operator $K$, using six ROIs.}
\label{tab:numerical_results_tlos}
\begin{tabular}{lcccc}
\hline
\textbf{Operator \( K \)} & \textbf{String} & \textbf{HD} & \( \gamma_3 \) & \( C_w \)  \\
\hline 
Clinical Data                   & TLOS & -- & -- & --  \\
Diffusion on   $\Gr_\textsc{c}$ & TLOS & 0 & 0.002 & 1.58  \\
Convolution on $\Gr_\textsc{l}$  & TSOL & 2 & 0.002 & 1.58   \\
Convolution on $\Gr_\textsc{nl}$ & TLOS & 0 & 0.002 & 1.58   \\
\hline
\end{tabular}
\end{table}
\begin{table}[h!]
\centering
\caption{Deterioration patterns and model parameters for different choices of the operator $K$, using ten ROIs.}
\label{tab:numerical_results_tfols}
\begin{tabular}{lcccc}
\hline
\textbf{Operator \( K \)} & \textbf{String} & \textbf{HD} & \( \gamma_3 \) & \( C_w \)  \\
\hline
Clinical Data & TFOLS & -- & -- & --\\
Diffusion on $\Gr_\textsc{c}$ & FTOLS & 2 & 0.001 & 50   \\
Convolution on $\Gr_\textsc{l}$  & TFOLS & 0 & 0.009 & 50  \\
Convolution on $\Gr_\textsc{nl}$  & FTOSL & 4 & 0.002 & 1.58   \\
\hline
\end{tabular}
\end{table}
From Tables~\ref{tab:numerical_results_tlos} and~\ref{tab:numerical_results_tfols}, we observe that the diffusion operator on the cumulative connectome $\Gr_\textsc{c}$ is able to reproduce only the deterioration pattern TLOS, corresponding to the configuration with six ROIs. Similarly, the convolution operator on the structural connectome $\Gr_\textsc{nl}$, where edge weights are defined as the ratio between the number of fibers and their length, also reproduces only the TLOS pattern.
On the other hand, the convolution operator on the length-based graph $\Gr_\textsc{l}$ is the only one capable of reproducing the deterioration pattern TFOLS, associated with the configuration involving ten ROIs.
Comparing the results obtained with 6 ROIs with those obtained with 10 ROIs, a sort of ``complementarity" stands out: what works best in one case works worse in the other. This shows that selecting the graph is an integral part of the modelling process and should be carried out with great attention.

To further highlight the importance of incorporating cumulative connectivity, we consider two alternative approaches: diffusion on the structural connectomes $\Gr_\textsc{nl}$ and $\Gr_\textsc{l}$, and convolution on $\Gr_\textsc{nl}$ and $\Gr_\textsc{l}$ with a kernel that does not exploit cumulative information. In the latter case,  the kernel $\bk^\text{sp}_\tau \in \mathbb{R}^N$ is constructed using shortest-path distances from each node to all others, and it is defined as:
\begin{equation}
    k^\textsc{sp}_\tau(i) = \sum_{j \in \Gr_\textsc{l}}  e^{\-\ell^\textsc{sp}_{ij} / \tilde{\delta}_\textsc{k}},\quad i=1,\ldots,N,
\end{equation}
where $\ell^\textsc{sp}_{ij}$ denotes the length of the shortest admissible path between vertices $i$ and $j$ (we consider only paths whose length is below the fixed threshold $R_\textsc{c}$) and $\tilde{\delta}_\textsc{k}=1$. This formulation does not account for the multiplicity of connections, and therefore lacks the cumulative nature that characterizes the current approach.
Tables~\ref{tab:numerical_results_tlos_mdpi} and~\ref{tab:numerical_results_tfols_mdpi} report the results obtained:
%using these approaches
none of these approaches successfully reproduces the deterioration pattern observed in the clinical data, 
neither in the case of the six ROIs nor in the case of the ten ROIs,
further supporting the relevance of incorporating cumulative connectivity into the model. 

%Nevertheless, none of these approaches successfully reproduces the deterioration pattern observed in the clinical data.

\begin{table}[h!]\ 
\centering
\caption{Deterioration patterns and model parameters for different choices of the operator $K$, using six ROIs.}
\label{tab:numerical_results_tlos_mdpi}
\begin{tabular}{lcccc}
\hline
\textbf{Operator \( K \)} & \textbf{String} & \textbf{HD} & \( \gamma_3 \) & \( C_w \)  \\
\hline
Clinical Data & TLOS & -- & -- & --  \\
Diffusion on $\Gr_\textsc{l}$ & LTOS & 2 & 0.001 & 0.5   \\
Diffusion on $\Gr_\textsc{nl}$ & LTOS & 2 & 0.001 &  0.5   \\
Convolution on $\Gr_\textsc{l}$   & TOSL &  3   & 50 & 0.5  \\
Convolution on $\Gr_\textsc{nl}$   & LSTO & 4 & 150   & 0.5  \\
\hline
\end{tabular}
\end{table}
\begin{table}[h!]
\centering
\caption{Deterioration patterns and model parameters for different choices of the operator $K$, using ten ROIs.}
\label{tab:numerical_results_tfols_mdpi}
\begin{tabular}{lcccc}
\hline
\textbf{Operator \( K\)} & \textbf{String} & \textbf{HD} & \( \gamma_3 \) & \( C_w \)  \\
\hline
Clinical Data                  & TFOLS & -- & -- & --  \\
Diffusion on $\Gr_\textsc{l}$  & FTOSL &  4 & 0.001  &  0.5    \\
Diffusion on $\Gr_\textsc{nl}$ & FTOLS &  2 & 0.001  &  0.5    \\
Convolution on $\Gr_\textsc{l}$  & FOLTS & 4 & 50 & 0.5  \\
Convolution on $\Gr_\textsc{nl}$  & FSTLO  & 4  & 150   & 0.5  \\
\hline
\end{tabular}
\end{table}
\section{Discussion and conclusion}
This study pursued two primary goals: first, to compare various modeling approaches for $\tau$ protein spread in the AD brain, specifically considering the synergistic presence of A$\beta$ protein; and second, to assess model validation by comparing generated numerical data with clinical data. Indeed, we firmly believe that a comparison between theoretical model results and clinical data is necessary to establish the practical utility of mathematical models in AD research.
To achieve this, we developed models where both proteins evolve on specialized networks derived from human connectomes and we introduced a new concept of deterioration pattern, enabling us to make the comparison mentioned above. The distinct physiological and biological characteristics of A$\beta$ and $\tau$ proteins necessitated different connectome structures. For A$\beta$, which propagates over short distances, we utilized a novel intrinsic proximity connectome with a standard diffusion Laplacian. For $\tau$ protein, hypothesized to spread prion-like over longer distances, we investigated distinct mathematical models across different networks: diffusion via Laplacian on a newly introduced cumulative connectome, and spreading through a convolution operator on two structural connectomes. An important novelty of this paper concerns the definition of a new convolution kernel, which takes into account the cumulative nature of the connections between brain regions and appears to provide the best performance in terms of comparison with clinical data. Indeed the numerical results obtained show that each operator $K$ is able to reproduce only one of the two clinical deterioration patterns, but not both, even when varying the model parameters.
In particular, the diffusion operator on the cumulative connectome $\Gr_\textsc{c}$ successfully reproduces the pattern TLOS, but fails to match TFOLS. Similarly, the convolution operator on the structural graph $\Gr_\textsc{nl}$ also reproduces only TLOS. On the other hand, the convolution operator on the length-based graph $\Gr_\textsc{l}$ 
is the only one capable of reproducing TFOLS. The fact that no operator can reproduce both patterns, even with different parameter settings, suggests that the mathematical structure of the model, including the underlying graph, must be carefully tailored to the specific phenomenon under investigation.

It is worth noting that in our previous work we also considered a convolution operator on the graph $\Gr_\textsc{nl}$, but with a different kernel. In that case, the model was not able to reproduce the TLOS pattern, further confirming that the choice of both the graph and the convolution kernel plays a crucial role in shaping the modeled dynamics.

As future work, we aim to apply the same comparative approach used for $\tau$ to validate the model's predictions for $A\beta$ dynamics against clinical data. This will allow us to further assess the model's ability to reproduce observed patterns of pathology and its consistency with medical observations.

\section*{Use of AI tools declaration}
The authors declare they have not used Artificial Intelligence (AI) tools in the creation of this article.

\section*{Institutional Review Board Statement} As per ADNI protocols, all procedures performed in studies involving human participants were in accordance with the ethical standards of the institutional and/or national research committee and with the 1964 Helsinki declaration and its later amendments or comparable ethical standards. More details can be found at adni.loni.usc.edu. (This article does not contain any studies with human 
participants performed by any of the authors)

\section*{Informed Consent Statement} Authors received the consent of publication from ADNI.

\section*{Data Availability Statement:} Publicly available datasets were analyzed in this study. This data can be found at these urls:
\url{http://adni.loni.usc.edu} and
\url{ https://braingraph.org } see \cite{szalkai2017parameterizable,szalkai2019high}.
\\
All data produced by the authors are available upon request from the authors.

\section*{Acknowledgments}
G.L. is member of the Gruppo Nazionale per il Calcolo Scientifico (GNCS) of the Istituto Nazionale di Alta Matematica (INdAM) and this work was partially supported by INdAM-GNCS under Progetti di Ricerca 2024 and 2025. 
M.C.T. is supported by PRIN 2022 F4F2LH - CUP J53D23003760006 ``Regularity problems in sub-Riemannian structures''.
A.S. and M.C.T. are members of the Gruppo Nazionale per l'Analisi Matematica, la Probabilit\`a e le loro Applicazioni (GNAMPA) of INdAM, Italy. 

Data collection and sharing for this project was funded by the Alzheimer's Disease Neuroimaging Initiative (ADNI) (National Institutes of Health Grant U01 AG024904) and DOD ADNI (Department of Defense award number W81XWH-12-2-0012). ADNI is fundedby the National Institute on Aging, the National Institute of Biomedical Imaging and Bioengineering, and through generous contributions from the following: AbbVie, Alzheimer’s Association; Alzheimer’s Drug Discovery Foundation; Araclon Biotech; BioClinica, Inc.; Biogen; Bristol-Myers Squibb Company; CereSpir, Inc.; Cogstate;
Eisai Inc.; Elan Pharmaceuticals, Inc.; Eli Lilly and Company; EuroImmun; F. Hoffmann-La Roche Ltd and its affiliated company Genentech, Inc.; Fujirebio; GE Healthcare; IXICO Ltd.; Janssen Alzheimer Immunotherapy Research \& Development, LLC.; Johnson \& Johnson Pharmaceutical Research \& Development LLC.; Lumosity; Lundbeck; Merck \& Co., Inc.; Meso Scale Diagnostics, LLC.; NeuroRx Research; Neurotrack Technologies; Novartis Pharmaceuticals Corporation; Pfizer Inc.; Piramal Imaging; Servier; Takeda Pharmaceutical Company; and Transition Therapeutics. The Canadian Institutes of Health Research is providing funds to support ADNI clinical sites in Canada. Private sector contributions are facilitated by the Foundation for the National Institutes of Health (www.fnih.org). The grantee organization is the Northern California Institute for Research and Education, and the study is coordinated by the Alzheimer’s Therapeutic Research Institute at the
University of Southern California. ADNI data are disseminated by the Laboratory for Neuro Imaging at the University of Southern California.

\section*{Conflict of interest}
The authors declare there is no conflict of interest.

%\bibliographystyle{plain}
%\bibliography{my_biblio}

\end{document}